\documentclass[11pt]{amsart}
\usepackage{amsmath, amsthm, latexsym, amssymb, amsfonts,
epsfig}

\newenvironment{pf}{\proof[\proofname]}{\endproof}
\theoremstyle{plain}
\newtheorem{Th}{Theorem}[section]
\newtheorem{Cor}[Th]{Corollary}

\newtheorem{Lemma}[Th]{Lemma}
\numberwithin{equation}{section} \theoremstyle{definition}

\newtheorem{Def}[Th]{Definition}
\newtheorem{Conj}[Th]{Conjecture}
\newtheorem{Remark}[Th]{Remark}

\newcommand{\cal}[1]{\mathcal{#1}}
\newcommand{\C}{\mathbb C}
\newcommand{\Z}{\mathbb Z}
\newcommand{\R}{\mathbb R}

\newcommand{\T}{{\mathbb T}}

\newcommand{\Vol}{\operatorname{Vol}}

\renewcommand{\Im}{\text{\rm Im }}
\renewcommand{\Re}{\text{\rm Re }}
\newcommand{\RP}{\mathbb RP}
\newcommand{\w}{\widetilde}
\renewcommand{\mod}{\mathrm{mod }}

\newcommand{\rl}[1]{Lemma~\ref{L:#1}}

\newcommand{\rr}[1]{Remark~\ref{R:#1}}

\newcommand{\re}[1]{(\ref{e:#1})}
\newcommand{\rcr}[1]{Corollary~\ref{C:#1}}
\newcommand{\rt}[1] {Theorem~\ref{T:#1}}

\newcommand{\rf}[1]{Figure~\ref{F:#1}}
\newcommand{\rconj}[1]{Conjecture~\ref{Conj:#1}}

\begin{document}


\title{Zeros of systems of exponential sums and trigonometric polynomials}
\author{Evgenia Soprunova}

\address{Department of Mathematics and Statistics,
University of Massachusetts, Amherst, MA 01003}
\email{soproun@math.umass.edu}

\begin{abstract} Gelfond and Khovanskii found a formula for the sum
of the values of a Laurent polynomial at the zeros of a system of
$n$ Laurent polynomials in $(\C^{\times})^n$ whose Newton polyhedra
have generic mutual positions. An exponential change of variables
gives a similar formula for exponential sums with rational
frequencies. We conjecture that this formula holds for exponential
sums with real frequencies. We give an integral formula which
proves the existence-part of the
conjectured formula not only in the complex situation but also in
a very general real setting. We also prove the conjectured formula
when it gives answer zero, which happens in most cases.
\end{abstract}
\maketitle
\section{Motivation and Summary}
Algebraic geometry is concerned with the study of zero sets of
algebraic polynomials $f(x)=\sum_\alpha c_\alpha x^\alpha$, where
$\alpha=(\alpha_1,\dots,\alpha_n)\in\Z^n$, $x=(x_1,\dots,x_n)$. If we allow the
exponents $\alpha$ to be real vectors, $f$ becomes a multi-valued
function. This can be remedied by an exponential change of
variables $x_i=\exp 2\pi z_i$, after which we obtain  a
single-valued exponential sum of the form $\sum_\alpha c_\alpha\exp 2\pi\alpha
z$, where $\alpha=(\alpha_1,\dots,\alpha_n)\in\R^n$ is the vector of frequencies, $z=(z_1,\dots,z_n)$
is the vector of variables, and $\alpha z$ is the standard scalar product.

Some results from algebraic geometry can be generalized to this
wider class of functions. For example,
Bernstein's theorem states
that the number of zeros of a generic system of $n$ algebraic
equations in $(\C^\times)^n$ with a fixed collection of Newton polyhedra is
equal to $n!$ times the mixed volume of the Newton polyhedra of
the system. A system of $n$ exponential equations in $n$ variables
usually has infinitely many isolated zeros, thus one has to study the distribution
of these zeros in order to obtain finite invariants.
O.~Gelfond proved that the mean number of complex zeros of a
system of $n$ exponential sums in $n$ variables with real
frequencies whose Newton polyhedra have sufficiently general
mutual positions (so-called developed collection of Newton polyhedra),
is equal to $n!$ times the mixed volume of the Newton
polyhedra of the system \cite{G2}.

O.~Gelfond and A.~Khovanskii found a formula for the sum
of the values of a Laurent polynomial at the zeros of a system of
$n$ Laurent polynomials in $\C^n$ with a developed collection of
Newton polyhedra \cite{GKh1,GKh2}. This formula splits into
two components. One of them is geometrical and reflects the mutual
positions of the Newton polyhedra of the system, while the other
component is expressed explicitly in terms of the coefficients of
the polynomial and the system.

An exponential change of variables gives a similar formula for
exponential sums with rational frequencies. We conjecture that
this formula also holds for exponential sums with real frequencies.
Here is some evidence for this. If the
exponential sum that we are summing up is identically equal to one, the
formula follows from combining two results: Gelfond's
generalization of Bernstein's theorem \cite{G2} and the new
formula for the mixed volume \cite{Kh1}. In \cite{So1} the
conjectured formula is proved in dimension one.

The conjectured formula, first of all, implies that the mean value
exists. We prove the existence not only in
the complex situation 
but also in a very general real setting by providing an integral
formula for the mean value.

If the frequencies of the exponential sum that we are summing up
are not commensurate with the frequencies of the system the
conjectured formula states that the mean value is equal to zero.
We show  that this is actually true and therefore prove the conjectured
formula in most cases. For example, if the exponential
sum that we are summing up is a single exponent $\exp 2\pi\alpha
z$, $z\in\C^n,\ \alpha\in\R^n$, then the formula is proved for all
values of $\alpha$ except for a countable set in $\R^n$.

Our arguments represent a combination of real analytic geometry
and ergodic theory. They are based on two theorems from completely
different parts of mathematics: the cell-decomposition theorem for
subanalytic sets (Appendix B) and Weyl's equidistribution law for
multidimensional trajectories in the real torus (Appendix A).
Weyl's equidistribution
law for one-dimensional trajectories is
 a classical theorem which was published in 1916 (see \cite{W3}).  The proof
of the corresponding law for multidimensional trajectories in \cite{So} is a
direct generalization of Weyl's original argument.

{\it Remarks.}
Similar ideas first appeared in Weyl's papers \cite{W1} and \cite{W2} where he solves the
mean motion problem.
S.~Gusein-Zade and A.~Esterov were dealing with related questions
and were using close techniques in \cite{GZ1,GZ2,E} to prove the
existence of the mean Euler characteristic and mean Betti numbers
of level sets and sets of smaller values of a quasiperiodic
function. Their motivation is completely different from ours: it
comes from the analysis of some models of chaotic behavior
appearing in quasicrystal structures. See also a recent paper \cite{Arnold}
of V.~I.~Arnold for a close discussion.

{\it Acknowledgements.}
This paper is a part of the author's Ph.D. thesis \cite{So} defended in spring 2002. I would
like to thank my thesis advisor Askold Khovanskii for stating
the problem and for his constant attention to this work.

\section{The Gelfond--Khovanskii formula and its conjectured generalization
to the case of exponential sums with real frequencies}

Let $\Delta_1,\dots,\Delta_n$ be convex polyhedra in $\R^n$ and
$\Delta=\Delta_1+\dots+\Delta_n$ be their Minkowski sum.
\begin{Def}
A collection of faces
$\{\Gamma_j:\Gamma_j\subset\Delta_j, j=1\dots n\}$ is {\it
coordinated} if there exists a nonzero linear function on $\R^n$
whose maximum on the polyhedron $\Delta_j$ is attained
exactly on the face $\Gamma_j\subset\Delta_j$ for all $j=1\dots n$.
A collection of polyhedra $\Delta_1,\dots,\Delta_n$ is
{\it developed} if none of the polyhedra is a vertex,
and in any coordinated collection of faces there is at least one vertex.
\end{Def}

Fix a system
$$P_1(z)=\dots=P_n(z)=0, \quad z\in\C^n$$ 
of $n$ Laurent equations with a developed collection of
Newton polyhedra $\Delta_1,\dots,\Delta_n$. The
Gelfond--Khovanskii formula states that the sum of the values of a
Laurent polynomial $Q$ over the zeros of the system in
$(\C^{\times})^n$ is equal to
$$(-1)^n\sum_{\alpha} k_\alpha C_\alpha,$$
where the summation is performed over the vertices $\alpha$ of
$\Delta=\Delta_1+\cdots+\Delta_n$, $k_{\alpha}$ is the
combinatorial coefficient that corresponds to the vertex $\alpha$
(the combinatorial coefficients are  integers that reflect the
mutual position of the Newton polyhedra
$\Delta_1,\dots,\Delta_n$ of the system, see \cite{GKh1,GKh2} for the definition), and 
$C_\alpha$ is an explicit Laurent polynomial in the coefficients of
$P=P_1\cdots P_n$ and $Q$. This result was announced in \cite{GKh1},
and a proof was given in \cite{GKh2}.

If we allow the exponents of the Laurent polynomials to be real vectors,
we obtain multi-valued functions. After an exponential
change of variables they become single-valued exponential sums
with real frequencies. Now we make some preparations before
we formulate
the conjectured generalization of the Gelfond--Khovanskii formula to the case
of exponential sums.

Let $\Lambda$ be a finite set in $\R^n$. An {\it exponential sum
with the spectrum} $\Lambda$ is a function $F:\C^n\to\C$ of the
form
$$F(z)=\sum_{\alpha\in\Lambda}c_{\alpha}\exp 2\pi\alpha z,
$$
where the summation is performed over the frequencies
$\alpha\in\Lambda$, $z=(z_1,\dots,z_n)\in\C^n$, $c_{\alpha}$
are nonzero complex numbers, and $\alpha z$ is the
standard scalar product. The {\it Newton polyhedron} of an
exponential sum is the convex hull $\Delta(\Lambda)$ of its
spectrum~$\Lambda$. 

Fix a system of exponential sums 
\begin{equation} \label{e:expsystem}
F_1(z)=\cdots=F_n(z)=0,\quad z\in\C^n
\end{equation}
with a developed collection of Newton polyhedra
$\Delta_1,\dots,\Delta_n$ in $\R^n$. Due to the assumption that the collection of the
polyhedra is developed, there exists $R>0$ such that all the zeros of the system lie in a strip
$S_R\times\text{Im }\C^n$, where $S_R\subset\text{Re }\C^n$ is a
ball of radius $R$ centered at the origin. This
implies that all the zeros are isolated (see \rt{Gelfond}). 

Let $G$ be another exponential sum with real frequencies.
The sum of the values at
the zeros of a system in the exponential case is replaced with the
result of averaging $G$ over the zeros of the system along
the imaginary subspace. Let $\Omega\subset\text{Im }\C^n$ be a bounded set with 
nonzero volume. For
$\lambda\in\R$, define $S_{\Omega}(\lambda)$ to be the sum of the
values of $G$ at the zeros of \re{expsystem} (counting multiplicities) that belong to
the strip $\R^n\times\lambda\Omega\subset\C^n$.

\begin{Def}
The {\it mean value} $M_{\Omega}$ of $G$ over the zeros
of the system \re{expsystem} is the limit of
$S_{\Omega}(\lambda)/\Vol(\lambda\Omega)$, as $\lambda$
approaches infinity.
\end{Def}

The Minkowski sum $\Delta=\Delta_1+\cdots+\Delta_n$ is the Newton
polyhedron of the product $F=F_1\cdots F_n$. The exponent
$\exp 2\pi\alpha z$ that corresponds to a vertex $\alpha$ in
$\Delta$ appears in $F$ with a nonzero
coefficient $d_\alpha$. Let $\w F=F/(d_\alpha \exp
2\pi\alpha z)$. The constant term of the exponential sum $\w F$
is equal to one. We define the exponential series for $1/\w F$ by
the formula 
$$1/\w F=1+(1-\w F)+(1-\w F)^2+\cdots.$$
 Since each exponent appears with a nonzero coefficient in a
finite number of terms, the coefficients of this series are
well-defined. Let $C_\alpha$ be the constant term in the
formal product of this series and $(1/d_\alpha)\exp(-2\pi\alpha z
)G\det(\frac{\partial F}{\partial z })$.

\begin{Conj}\label{Conj:GKhexp}
The mean value $M_{\Omega}$
of an exponential sum $G$ over the zeros of the system \re{expsystem} is
equal to 
$$\frac{1}{(-2\pi)^n}\sum_{\alpha} k_\alpha C_\alpha,$$ where
the summation is performed over the vertices $\alpha$ of the
Minkowski sum $\Delta=\Delta_1+\cdots+\Delta_n$, $C_\alpha$ are the
constant terms of the series as defined above, and $k_\alpha$ are the
combinatorial coefficients.   
\end{Conj}
As we mentioned in the summary, this formula is proved in the following 
three cases:
\begin{itemize}
\item[1.]The frequencies $\alpha$ are rational (the formula is obtained
by an exponential change of variables from the algebraic case).
\item[2.]$G=1$ (the formula follows from combining two
results: Gelfond's generalization of Bernstein's theorem
\cite{G2} and the new formula for mixed volume \cite{Kh1}).
\item[3.]$n=1$ (see \cite{So1}).
\end{itemize}

In the exponential case even the existence of the mean value is not 
obvious. We give an integral formula for the mean value which proves 
the existence not only in the complex situation
(\rt{compmain}) but also in a very general real setting
(\rt{main}).

Assume that the frequencies of the exponential sum that we are
summing up are not commensurate with the frequencies of the
system. Then there is no constant term in the series defined above, that is, all
$C_{\alpha}$ are equal to zero, and the
conjectured formula states that the mean value is equal to zero.
Using the integral representation, we show in \rcr{nolik} that this is actually true. This proves the
conjectured formula in most cases. For example, if the exponential
sum that we are summing up is a single exponent $\exp 2\pi\alpha
z$, then we proved the formula for all values of $\alpha$ except
for a countable set in $\R^n$.

\section{Averaging over the isolated
intersections of a subanalytic set with a dense orbit on $\T^N$}

Here we present a construction and formulate a theorem on which all our
results are based. Let $\w V$ be a subanalytic subset of the real torus
$\T^N$ and $\cal O$ a dense multidimensional orbit on $\T^N$.
We prove an integral formula for the mean value of a bounded subanalytic function on
$\w V$ over the isolated intersections of the set $\w V$ and the
orbit $\cal O$. This formula, in particular, implies that the mean value
always exists.

Let $x=(x_1,\dots,x_n)\in\R^n$ and
$\varphi=(\varphi_1,\dots,\varphi_N)\in\T^N=\R^N/\Z^N$. A
linear map $\Phi:\R^n\to\R^N$ defines an action of $\R^n$ on $\T^N$
by $$x:\varphi\mapsto(\varphi+\Phi(x))\,\mod\,\Z^N.$$
Let $\cal O$ be an orbit of this action. We require that the orbits are
dense in the torus, which means that there are
no integral vectors orthogonal to the $n$-plane
$\Phi(\R^n)\subset\R^N$ (Appendix A).
We will also assume that $\Phi$ is injective, that is, the orbits are
$n$-dimensional.

Let $\w V$ be a subanalytic subset of $\T^N$ and $\w T$ a
bounded subanalytic function on $\w V$. Set $V:=\Phi^{-1}(\w
V\cap\cal O)\subset\R^n$, and define a function $T(x)=\w
T(\Phi(x))$ on $V$.
Let $\Omega$ be a bounded subset of $\R^n$ with 
nonzero volume. For $\lambda>0$,
define $S_{\Omega}(\lambda)$ to be the sum of
the values of $T$ at the isolated points of $V$ that belong to
$\lambda\Omega$.

\begin{Def}
The {\it mean value} $M_{\Omega}$ of $\w T$ over the
isolated  intersections  of $\w V$ and
$\cal O$ (in the topology of $\cal O$) is the limit:
$$M_{\Omega}=\lim_{\lambda\to\infty}
\frac{S_{\Omega}(\lambda)}{\Vol(\lambda\Omega)}.$$
\end{Def}

Let $p:\R^N\to L$ be the orthogonal projection to the linear
subspace $L\subset\R^N$ orthogonal to $\Phi(\R^n)$. Define
$M_{N-n}(\w V)$ to be the smooth $(N-n)$-dimensional part of $\w V$,
transversal to the orbit. Here is the precise three step definition
(see Figure 1 below where we assume that the direction of the orbit is vertical):

\begin{enumerate}
\item  Let $\w V_0$ be the set of all points $\varphi\in\w V$
such that the intersection of $\w V$ with the $n$-plane through
$\varphi$ parallel to the orbit is locally just the point
$\varphi$ itself.
\item  Let $\w V_1$ be the set of all points $\varphi\in\w
V_0$ such that $\w V_0$ is a $C^1$-manifold of dimension $N-n$ in
a small neighborhood of $\varphi$.
\item  Let $M_{N-n}(\w V)$ be the set of all points $\varphi\in\w
V_1$ such that the projection $p$ is regular at each of these
points.
\end{enumerate}

\begin{figure}[h]
\centerline{\epsfysize = 1.35 in \epsffile{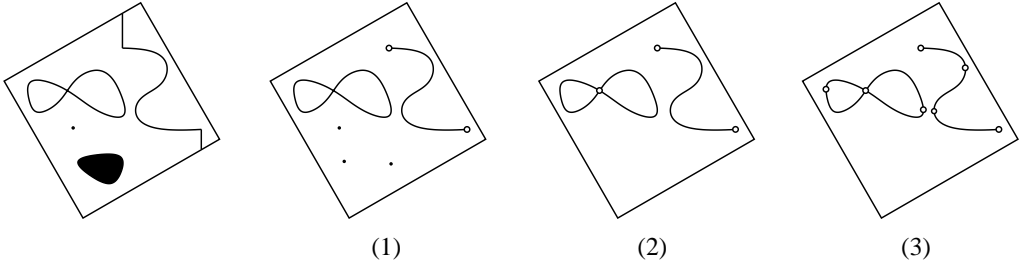}}
\caption{Definition of $M_{N-n}(\w V)$}
\end{figure}

Let $A^j=\Phi(e_j)$, where $e_j$ is the $j$-th vector in the 
standard basis for $\R^n$. The quotient of the standard volume form
on $\R^N$ by the image under $\Phi$ of the standard volume form on
$\R^n$ defines a volume form
$\omega=d\varphi(\dots,A^1,\dots,A^n)$ and an orientation on $L$.
Then $M_{N-n}(\w V)$ is a $C^1$-manifold, and the projection $p$ defines an
($N-n$)-form $p^*\omega$ on
$M_{N-n}(\w V)$.

\begin{Def}
The {\it $\cal O$-transversal $\w T$-weighted
volume} of $\w V$ is the integral of the form $\w
Tp^*\omega$ over $M_{N-n}(\w V).$
\end{Def}

\begin{Def}
We say that the {\it fractal dimension} of a set $X$ is less than
or equal to $n$ if there exist constants $c$ and $\epsilon_0>0$
such that for any $0<\epsilon<\epsilon_0$, the set $X$ can be
covered by no more than $c\epsilon^{-n}$ balls of radius
$\epsilon$. We say that the fractal dimension is less than $n$ if
it is less than or equal to $n-\delta$ for some positive $\delta$.
\end{Def}
Our main result is the following theorem:
\begin{Th}\label{T:general}
Let $\cal O$ be an $n$-dimensional dense
orbit in the torus $\T^N=\R^N/\Z^N$, $n<N$, $\w V$ a
subanalytic subset of $\T^N$, and $\w T$  a bounded
subanalytic function on $\w V$. Let $\Omega$ be a bounded subset
of $\R^n$ with
nonzero volume whose boundary has fractal dimension less than $n$.

Then the mean value $M_{\Omega}$ of $\w T$ over the
isolated intersections  of $\w V$ and
$\cal O$ always exists and equals the $\cal O$-transversal
$\w T$-weighted volume of $\w V$.
\end{Th}

We prove this theorem in Section 6.

\section{Zeros of systems of trigonometric polynomials}
We apply \rt{general} to prove that the mean value of a
quasiperiodic trigonometric polynomial $T$ over the isolated
points of a quasiperiodic semitrigonometric set $V$ (a set
described by quasiperiodic trigonometric polynomials, see the
precise definition below) always exists. We obtain an integral
representation for the mean value from which we conclude that in the case when the
frequencies of $T$ are not commensurate with the frequencies of
the quasiperiodic trigonometric polynomials that describe $V$, the
mean value is equal to zero.

\subsection{Definitions}
Let  $x\in\R^n$. A {\it quasiperiodic trigonometric
polynomial} is a trigonometric polynomial of the form
$$T(x)=\sum_{k=1}^{p}c_k\cos 2\pi\alpha_kx+d_k\sin 2\pi\alpha_kx$$
with real coefficients $c_k, d_k\in\R$ and real frequencies
$\alpha_k\in\R^n$. Here $\alpha_kx$ is the standard scalar product.

A subset of $\R^n$ is called quasiperiodic semitrigonometric if it is described
by quasiperiodic trigonometric equations and inequalities. More
precisely, we say that $V\subset\R^n$ is {\it quasiperiodic
semitrigonometric} if it
can be represented in the form:
$$V
=\bigcup_{i=1}^s\bigcap_{j=1}^{r_i}V_{ij},$$ where each $V_{ij}$
is either $\{x\in\R^n\,|\,T_{ij}(x)=0 \}$ or $\{x\in\R^n\,|\,T_{ij}(x)>0
\}$, for some quasiperiodic trigonometric
polynomials $T_{ij}$.

Let $T$ be
a quasiperiodic trigonometric polynomial, $V$ a quasiperiodic semitrigonometric set in $\R^n$ described by $T_{ij}$,  
and $\Omega$ a bounded subset of
$\R^n$ with nonzero volume. For $\lambda>0$, define $S_{\Omega}(\lambda)$
to be the sum of the values of $T$ at the isolated points
of $V$ that belong to $\lambda\Omega$. The {\it mean value} $M_{\Omega}$ of $T$ over the
isolated points of $V$ is defined by
$$M_{\Omega}=\lim_{\lambda\to\infty}\frac{S_{\Omega}(\lambda)}{\Vol(\lambda\Omega)}.$$

\subsection{Averaging a quasiperiodic trigonometric polynomial over
the isolated points of a quasiperiodic semitrigonometric set}

Consider the set of all frequencies of $T$ and 
$T_{ij}$. This set
generates a subgroup $\cal A$ of $(\R^n,+)$, which is a direct
sum of infinite cyclic subgroups:
$$\cal A = (A_1)\oplus\cdots\oplus(A_N).$$
That is, there exist $A_1,\dots,A_N\in\R^n$ with no integral relations,
such that each of the frequencies is an integral combination of
$A_1,\dots,A_N$.

Define a linear map $\Phi:\R^n\to\R^N$ by $\Phi(x)=(A_1x,\dots,A_Nx)$.
Then $\R^n$ acts on $\T^N=\R^N/\Z^N$ by
$$x:\varphi\mapsto(\varphi+\Phi(x))\,\mod\,\Z^N.$$
Let $\cal O$ be the orbit of this action through the origin.
Since each of the frequencies $\alpha_k$ is an integral combination
of $A_1,\dots,A_N$, the functions $T$ and $T_{ij}$
are the restrictions to $\cal O$ of some
trigonometric polynomials $\w T$ and $\w T_{ij}$ defined
on $\T^N$.
For example, $\cos (2\pi\alpha x)$ is the restriction of $\cos 2\pi\sum m_i\varphi_i$,
where $\alpha=\sum m_iA_i$ for $m_i\in\Z$.

To apply \rt{general} here we need the
orbit to be dense (\rl{nointegral}), the mapping $\Phi$ to be injective
(\rl{lessthann}), and the dimension $n$ of the
orbit to be less than $N$ (if $\Phi$ is injective
then $n\leq N$; $N=n$ corresponds to the periodic case which we treat in \rr{periodic}).

\begin{Lemma}\label{L:nointegral}
There are no nonzero integral vectors orthogonal to the plane of
the orbit $\cal O$. Therefore, $\cal O$ is
dense in the torus $\T^N$.
\end{Lemma}
\begin{pf}

Let $A$ be the $N\times n$ matrix whose rows are the generators
$A_1,\dots,A_N$ of $\cal A$. We denote the columns of this matrix
by $A^1,\dots,A^n$. These columns generate the $n$-plane $\Phi(\R^n)$
of the orbit $\cal O$. If there exists a nonzero integral
vector $(k_1,\dots,k_N)$
orthogonal to $\cal O$ then it is orthogonal to
the vectors $A^1,\dots,A^n$.
This implies that $k_1A_1+\cdots+k_NA_N=0$, i.e. there exists
a nontrivial linear combination with integral coefficients on the
generators $A_1,\dots,A_N$ of $\cal A$, which gives a
contradiction. By Appendix A the orbit $\cal O$ is dense in
the torus $\T^N$.
\end{pf}

\begin{Lemma}\label{L:lessthann}
If the dimension of the orbit $\cal O$ is less than $n$, then
the set $V$ described by the $T_{ij}$ has no isolated points,
and the mean value $M_\Omega$ is equal to zero.
\end{Lemma}
\begin{pf}
If the dimension of the orbit is less than
$n$, then $A^1,\dots,A^n$ are linearly dependent, and
the solution space of the system $A^1x_1+\dots+A^nx_n=0$ is nontrivial.
If $T_{ij}(x)=0$, then
$T_{ij}(x+y)=0$ for each $y$ in the solution space. This means that $V$ has no
isolated points, and the mean value is equal to zero.
\end{pf}

If the orbit is $n$dimensional and $n<N$,
then $\R^n$ is mapped bijectively to a dense
$n$-dimensional orbit on the real torus $\T^N$ and trigonometric polynomials
$T_{ij}$ regarded as functions $\w T_{ij}$ on $\T^N$ define a set
$\w{V}$. The intersection of the set $\w{V}$ and the orbit
$\cal O$ is the semitrigonometric set $V$. We are adding up
the values of $\w T$ over the isolated (in the topology
of the orbit) points of intersection of the set $\w{V}$ and
the orbit. Here is the main result of this chapter which
now follows directly from \rt{general}.

\begin{Th}\label{T:main}
Let $V$ be a quasiperiodic semitrigonometric subset of $\R^n$,
$T$ a quasiperiodic trigonometric polynomial. Let
$\w V$, $\w T$, $\cal O$ be the corresponding set, function on
$\T^N$, and orbit, as constructed above. Let  $\Omega$ be a
bounded subset of $\R^n$ with nonzero volume, whose boundary has fractal dimension
less than $n$.

Then the mean value $M_\Omega$ always exists and is equal to the
$\cal O$-transversal\linebreak $\w T$-weighted volume of $\w V$,
assuming that the dimension of the orbit is $n$, and $n<N$. If the dimension
of the orbit is less than $n$, the mean value is equal to
zero.
\end{Th}
\begin{Remark}\label{R:periodic}
If the orbit has the same dimension as the torus, i.e.,
$N=n$, the set $V\subset\R^n$ and the function $T$
have $n$ linearly independent periods $A^1,\dots,A^n$, which
define a torus $\R^n/(\Z A^1\oplus\cdots\oplus\Z A^n)$. The mean
value $M_\Omega$ in this case is equal to the sum of the values of
$T$ at the isolated points of $V$ in this torus,
divided by the volume of the torus.
\end{Remark}
\begin{Remark}
This theorem holds and our proof works if $V$ is a quasiperiodic subanalytic
subset of $\R^n$, that is, V is the intersection of a subanalytic subset
$\w V$ of $\T^N$ with some dense $n$-dimensional orbit, and $T$ is
the restriction to this orbit of some bounded subanalytic function on
$\w V$.
\end{Remark}

\begin{Cor}\label{C:trigzero}
If all the multiples $k\alpha$, $k\in\Z\setminus\{0\}$ of each frequency
$\alpha$ of $T$ do not belong to the subgroup $\cal A$ of $(\R^n,+)$
generated by the frequencies of $T_{ij}$,
then the mean value $M_{\Omega}$ is equal to zero.
\end{Cor}
\begin{pf}
It is enough to check this statement for $T$ with a single
frequency $\alpha$, i.e., for $T(x)=c\cos 2\pi\alpha x+d\sin 2\pi\alpha x$. The set of
all frequencies of the trigonometric polynomials $T_{ij}$
describing the set $V$ generates a group $\cal A'$ in $(\R^n,+)$,
which is the direct sum of infinite cyclic subgroups:
$$\cal A '= (A_1)\oplus\cdots\oplus(A_{N-1}).$$
Since no multiple of the frequency $\alpha$ of $T$ belongs to $\cal
A'$, the group $\cal A$ obtained from $\cal A'$ by throwing in
$\alpha$,  is the direct sum:
$$\cal A = (A_1)\oplus\cdots\oplus(A_{N-1})\oplus(\alpha).$$

Therefore, raising the situation to the torus $\R^N/\Z^N$, we
obtain a function $\w T(\varphi_N)=c\cos \varphi_N+d\sin
\varphi_N$ that depends only on the variable $\varphi_N$, and a
set $\w V$, which is described by the functions $\w T_{ij}$ that do not depend on $\varphi_N$.

By \rt{main}, the mean value is equal to the $\cal O$-transversal $\w T$-weighted 
 volume of $\w V$. Since the $\w T_{ij}$ 
do not depend on $\varphi_N$, the manifold
$M_{N-n}(\w V)$ (which is the smooth $N-n$ dimensional,
transversal to the orbit $\cal O$, part of $\w V$) is the
direct product $\w V'\times S^1$, where $\w V'$ is a manifold in
$\R^{N-1}/\Z^{N-1}$, and $S^1=\R/\Z$. The form $\omega=p^\ast
d\varphi_1\dots d\varphi_N(\dots,A^1,\dots,A^n)$ can be written as
$\omega_1\wedge d\varphi_N$ where $\omega_1$ is a form on $\w V'$.
Therefore, the mean value is equal to
$$
\int_{M_{N-n}(\w V)} T\omega=\int_{\w V'}\left(\int_{S^1}(c\cos
2\pi\varphi_N+d\sin 2\pi\varphi_N) d\varphi_N\right)\omega_1=0.
$$
\end{pf}
\section{Zeros of systems of exponential sums}
In this section we apply our main result (\rt{general}) to
the complex case, that is, to computing the mean value of an
exponential sum over the zeros of a system of $n$ exponential sums
in $\C^n$. We show that this mean value always exists. From the
integral representation for the mean value that we obtain, we
deduce a proof of \rconj{GKhexp} in the
case when the frequencies of the exponential sum are not
commensurate with the frequencies of the system.

\subsection{Definitions}
Fix a system
\begin{equation}\label{e:system}
F_1(z)=\cdots=F_n(z)=0,\quad z\in\C^n
\end{equation}
of exponential
sums with a developed collection of Newton polyhedra. Let $\cal F$
be a family of systems of $n$ exponential sums with the same collection
of spectra as \re{system}, and whose coefficients have
the same absolute values as the corresponding coefficients in \re{system}.
The following theorem is a particular case of a result proved by Gelfond
\cite{G1,G2}. 

\begin{Th}\label{T:Gelfond}
There exists $R>0$ such that all the zeros of the
systems from the family $\cal F$ belong to the strip
$S_R\times\rm{Im}\,\C^n$, where $S_R\subset\rm{Re}\,\C^n$ is a ball
of radius $R$ centered at the origin. This implies
that all the zeros of the systems from $\cal F$ are
isolated.
\end{Th}

Let $G$ be an exponential sum with real frequencies, and $\Omega$  a bounded
subset of $\text{Im}\,\C^n$ with
nonzero volume. For $\lambda>0$, let $S_{\Omega}(\lambda)$
be the sum of the values of  $G$ at the zeros of \re{system}
(counting multiplicities) that belong to
the strip $\R^n\times\lambda\Omega\subset\C^n.$ Let
$S_{\Omega}(\lambda)^{\geq k}$ be the
sum of the values of  $G$ at the zeros (not counting multiplicities)
of the system \re{system} of multiplicity at least $k$ that belong to
the strip $\R^n\times\lambda\Omega\subset\C^n$.
By \rt{Gelfond},
$S_{\Omega}(\lambda)$ and $S_{\Omega}(\lambda)^{\geq k}$ are well--defined.

Define $$M_{\Omega}=\lim_{\lambda\to\infty}\frac{S_{\Omega}(\lambda)}{\Vol(\lambda\Omega)}, \ \ \
M_{\Omega}^{\geq k}=\lim_{\lambda\to\infty}\frac{S_{\Omega}^{\geq k}(\lambda)}{\Vol(\lambda\Omega)}.$$

\rt{compmain} below provides an integral formula for computing
the mean values $M_{\Omega}$, $M_{\Omega}^{\geq k}$ similar to the
formula in the real case (\rt{main}). In particular, \rt{compmain}
implies that the mean
value $M_{\Omega}$ always exists.
Another consequence  of this integral formula is
that in the case when the
frequencies of $G$ are not commensurate with the frequencies
of the exponential polynomials from the system, the
mean value $M_{\Omega}$ is equal to zero (\rcr{nolik}).
\subsection{Averaging an exponential sum over the zeros of a system of exponential sums}
As in the real case, we deduce an integral formula for the
mean value from \rt{general}. The exponential sums $F_j$, $G$
are quasiperiodic along the imaginary subspace of $\C^n$. Since the system
\re{system} has a developed collection of Newton polyhedra, the real parts of the zeros of this
system belong to some ball $S_R\subset\text{Re }\C^n$. Therefore, we can
restrict our attention to this ball, and think of the $F_j$ and $G$
as of functions quasiperiodic along the imaginary subspace and periodic
along the real subspace. This allows us to consider $\C^n$ as an orbit
on some real torus and apply \rt{general} to this situation.

As before, the set of all frequencies $\alpha$ of the
exponential sums $F_{j}$, $G$ generates a subgroup
$\cal A$ in $(\R^n,+)$, which is a direct sum of infinite cyclic subgroups: 
$$\cal A = (\alpha_1)\oplus\cdots\oplus(\alpha_N).$$
In other words, there exist $\alpha_1,\dots,\alpha_N\in\R^n$ with no integral
relations, such that each of the frequencies $\alpha\in\cal A$ is
an integral combination of $\alpha_1,\dots,\alpha_N$.

Define a linear map 
$$\Phi_1:\text{Im}\,\C^n\to\R^N\quad
\text{by}\quad\Phi_1(y)=(\alpha_1y,\dots,\alpha_Ny),\quad y\in\!\text{Im}\,\C^n.$$
According to \rt{Gelfond} the real
parts of the zeros of the system
\re{system} belong to the open ball
$S_R\subset\text{Re}\,\C^n$ centered at the origin.
Let $\Phi_0\!:\text{Re}\,\C^n\to\R^n$ be a composition of a shift and a rescaling
that maps $S_R$ inside the unit cube. We now think of $F_j$ and $G$
as of functions periodic along the real subspace.

Then $\Phi=\Phi_0\times\Phi_1$ maps
$\C^n=\Re\C^n\times\Im\C^n$ to $\R^n\times\R^N$. The map
$\Phi$ defines an action of $\C^n$ on $\T^{n+N}$ by
$$\varphi\mapsto(\varphi+\Phi(x,y))\,\mod\,\Z^{n+N},\ \ \varphi\in\T^{n+N},\ \ x\in\Re\C^n,\ \ y\in\Im\C^n.$$ 
Let $\cal O$ be the orbit of this action through the origin. Since each of the frequencies $\alpha$
is an integral combination  of $\alpha_1,\dots,\alpha_n$, the sums $F_j$, $G$
are the restrictions to $\cal O$ of some functions
$\w F_j$, $\w G$ defined on $\T^{n+N}$. We note that these
functions may only be non-analytic at the points $\varphi$ where at
least one of the coordinates $\varphi_1,\dots,\varphi_n$
is equal to zero.

The set of all zeros of the system \re{system} of multiplicity
at least $k$ is an analytic set defined by
exponential equations $F_j(z)=0$ and a few more exponential
equations whose frequencies belong to the subgroup $\cal A$ of
$(\R^n,+)$ (see, for example,
\cite{GKh}).

Let $V^{\geq k}$ be the set of all zeros of the system of
multiplicity greater than or equal to $k$. Consider the system of
exponential sums that defines $V^{\geq k}$. Regarded as functions on $\T^{n+N}$,
these sums define a subset $\w{V}^{\geq k}$ of the torus $\T^{n+N}$.
We now state the main result of this section.

\begin{Th}\label{T:compmain}
Consider a system of exponential sums with real frequencies
\begin{equation}\label{e:compmainsys}
F_1(z)=\cdots=F_n(z)=0,\quad z\in\C^n 
\end{equation}
with a developed collection of Newton polyhedra. Let $G$ be an exponential
sum with real frequencies, and  $\Omega$ a bounded subset of
$\R^n$ with
nonzero volume whose boundary has fractal dimension less than $n$.
Let $\w{V}^{\geq
k}$, $\w G$, $\cal O$ be the corresponding
subsets of $\T^{n+N},$ function on $\T^{n+N}$,
and orbit, as constructed above.

Then if $n<N$ the mean value $M_\Omega^{\geq k}$ is equal to the
$\w G$-weighted $\cal O$-transversal volume of
$\w{V}^{\geq k}$. The mean value
$M_\Omega$ is equal to the sum $M_{\Omega}^{\geq 1}+\dots+M_{\Omega}^{\geq k}$ for some~$k$.
\end{Th}
\begin{pf}
We first show that the dimension of the orbit $\cal O$ is $2n$.
It is enough to show that the dimension of the plane
$\Phi_1(\text{Im }\C^n)$ is $n$. This plane is generated by the
columns of the matrix of the linear map $\Phi_1$. The rows of this matrix
are the generators $\alpha_1,\dots,\alpha_N$ of the group $\cal
A$. Thus the dimension of the plane $\Phi_1(\text{Im }\C^n)$ is
equal to the dimension of the linear space generated by the
frequencies $\alpha_1,\dots,\alpha_N$. Since the collection of
the Newton polyhedra of the system \re{compmainsys} is
developed, its frequencies generate $\R^n$. We conclude that the dimension of the
orbit is $n$.

Similarly to the real case, it is easy to see that
there are no integral vectors orthogonal to the orbit $\cal
O$. Therefore, the orbit $\cal O$ is dense in the torus
$\T^{n+N}$ (see \rl{nointegral}).

To derive \rt{compmain} from \rt{general} we show that $\w
{V}^{\geq k}$ is an analytic subset of $\T^{n+N}$ and $\w
G$ is analytic on $\w{V}^{\geq k}$.
Each of $\w{V}^{\geq k}$ is defined by a system of equations
analytic everywhere on $\T^{n+N}$ but possibly at the points where at least
one of the coordinates $\varphi_1,\dots,\varphi_n$ is equal to
zero. Notice that $\w{V}^{\geq k}$ is contained in the
``cylinder'' $\Phi_0(S_R)\times\T^N$. Indeed, the intersections of
$\w V$ with the orbit through the origin belong  to this cylinder
since the zeros of the initial system lie
in the strip $S_R\times\text{Im }\C^n$. Restricting the functions
$\w F_j$ to a shifted orbit we obtain functions that differ from
$F_j$ by a shift along the imaginary subspace. By \rt{Gelfond}
the zeros of this new system belong to the strip
$S_R\times\text{Im }\C^n$, and we conclude that $\w{ V}^{\geq k
}\subset\Phi_1(S_R)\times\T^N$. Therefore, $\w{V}^{\geq k}$ is
locally defined by a system of analytic equations, and $G$
is analytic on $\w V^{\geq k}$. 

To show that $M_\Omega$ can be computed by the formula
$M_\Omega=M_{\Omega}^{\geq 1}+\dots+M_{\Omega}^{\geq k}$ for some $k$,
we need to prove that the sets $V^{\geq k}$ are empty, starting from some $k$.
Assume that there exists a sequence $z_k=x_k+iy_k$, $x_k\in
S_R$, such that $z_k$ belongs to $V^{\geq k}$. Raising this
sequence to $S_R\times\T^N$, we obtain a sequence $\w z_k$ in
$S_R\times\T^N$  such that $\w z_k$ belongs to $\w{V}^{\geq
k}$. Since $S_R\times\T^N$  is compact, there exists a limit point
$\w z_0\in S_R\times\T^N$. This point $\w z_0$ belongs to
each $\w{V}^{\geq k}$. Although this point does not necessarily
belong to the orbit $S_R\times\Phi_1(\text{Im }\C^N/\Z^N)$, it
belongs to some shifted orbit. The restrictions of the functions
$\w {F}_j$ to this shifted orbit define an exponential
system whose zeros are isolated and belong to $S_R\times\text{Im }{\C^n}$ 
(\rt{Gelfond}). Therefore, $\w z_0$ is a zero of the
shifted system of infinite multiplicity, which is impossible since
all the zeros of this system are isolated. We proved that all
$V^{\geq k}$ are empty starting from some $k$.

Now the theorem follows from \rt{general}.
\end{pf}

The proof of the following corollary is a repetition of the proof of \rcr{trigzero}.

\begin{Cor}\label{C:nolik}
If all the multiples $k\alpha,\ k\in\Z\setminus\{0\}$ of each
frequency $\alpha$ of $F$ do not belong to the subgroup $\cal A$
of $(\R^n,+)$ generated by the frequencies of the $F_j$, then the mean value $M_{\Omega}$ is
equal to zero.
\end{Cor}

\section{Proof of \rt{general}}
\begin{pf}Recall that $\cal O=(a+\Phi(\R^n))\,\mod\,\Z^N$, where $a\in\T^N$, and
$\Phi:\R^n\to\R^N$ is a linear injective map. We send the torus
$\T^N$ to the cube $I=[0,1)^N\subset\R^N$ by cutting $\T^N$ along
$\varphi_j=1$. Choose a new basis
$\{B^1,\dots,B^{N-n},A^1,\dots,A^n\}$ in $\R^N$ whose last $n$
coordinate vectors are the images of the standard basis in $\R^n$ under $\Phi$,
$A^j=\Phi(e_j)$. Note that the vectors $A^1,\dots,A^n$ generate the plane of
the orbit $\Phi(\R^n)$. Let the remaining vectors $B^1,\dots,B^{N-n}$ 
form a basis of the
$(N-n)$-plane $L$ orthogonal to the orbit, so that
the determinant of the change of coordinates matrix from the standard
basis in $\R^N$ to this new basis is equal to one.

Let $x=(x_1,\dots,x_n)$ be the coordinates along the orbit,
and $y=(y_1,\dots,y_{N-n})$ the coordinates in the orthogonal
plane $L$. We choose a cube $Q=Q_y\times Q_x$ whose edges go
along the new coordinate axes and $I\subset Q$. We now regard
$\w V$ as a subset of $I\subset Q$ and $\w T(\varphi)$
as a function $\w T(y,x)$ on $\w V\subset I$.

The set $\w V$ is a subset of
$Q$, subanalytic in $\RP^N$, and the graph of the function $\w
T$ is a subset of $Q\times\R $ subanalytic in $\RP^{N+1}$.
Therefore, by the cell decomposition theorem
(\rt{celldecomposition}), there exists a decomposition of
$\R^N$ with the fixed order of coordinates $(y,x)$ into finitely many cells such that $\w V$ is
a union of cells and $\w T$ is $C^1$ on each of these cells.

Since the cell decomposition respects the standard projection,
each cell is either a part of the closure of another cell, or is
disjoint from it. Only 
$(i_1,\dots,i_{N-n},0,\dots ,0)$-cells that do not belong to the
closure of a cell of higher dimension may have isolated
intersections with the orbit. We will show later that cells of
dimension less than $N-n$ do not contribute to the mean value.
Therefore, we can restrict our attention to 
$(1,\dots,1,0,\dots,0)$-cells with ${N-n}$ ones and $n$ zeros that do
not belong to a closure of a cell of higher dimension.

Let $C$ be one of such cells. Then
$$C =\{(y,x)\ |\
y\in D,\ x=H(y) \},$$ where $H=(H_1,\dots,H_n)$ for some $C^1$
functions $H_1,\dots,H_n$, and $D$ is a $(1,\dots,1)$-cell  in
$Q_y$.

Now we define a counting function $g$ on the torus
$\T^N$. Through each point $a\in C$ we draw an $n$-plane parallel
to the orbit. If $\varphi$ belongs to $(a+\Phi(B_r))/\Z^N$, where
$B_r\subset\R^n$ is the ball of radius $r$ centered at the origin,
we define $g(\varphi)=\w T(a)/v_r$, where $v_r$ is the volume of
the ball $B_r\subset\R^n$. (Here $r$ is chosen so small that the
mapping $x\mapsto\Phi(x)/\Z^N$ restricted to $B_r$ is one-to-one.)
Otherwise we define $g(\varphi)=0$.

In $(y,x)$-coordinates, $g$ is defined by
$$g(y,H(y)+x)=
\begin{cases}
    \w T(y,H(y))/v_r,& \text{if } y\in D,\ x\in B_r  \\
    0, & \text{otherwise}.
  \end{cases}
$$
\rf{counting} illustrates the definition of $g$ in the case when the cell
$C$ is one-dimensional. Through each $a\in C$ draw a vertical segment of length
$2r$ centered at $a$. This segment is the ball $a+B_r$ from the general construction.
Define $g$ to be identically equal to $\w T(a)/2r$ on that segment.

\begin{figure}[h]
\centerline{
\scalebox{0.8}{
\input{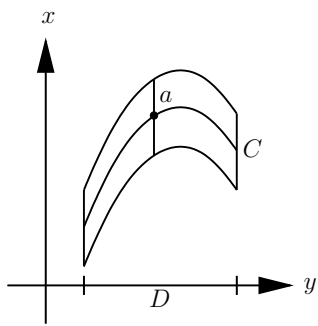}
             }
           }
         \caption{Definition of the counting function $g$}
         \label{F:counting}
       \end{figure}
\setlength{\unitlength}{1mm}

 Since the function ${\w T}$ is
continuous and bounded on $C$ and the functions defining $C$
are continuously differentiable, $g$ is Riemann
integrable on the torus $\T^N$. We apply Weyl's equidistribution
law for multidimensional trajectories (\rt{Weyl}):
\begin{equation}\label{e:weyllaw}
\lim_{\lambda\to\infty}\frac{1}{\Vol(\lambda\Omega)}
\int_{\lambda\Omega}g(\Phi(x))dx =\int_{\T^N}g(\varphi)d\varphi.
\end{equation}

Let us show first that the left hand side part of this identity is
exactly the contribution of the cell $C$ to the mean value
$M_{\Omega}$. If $a+B_r$, a ball centered at the point $a\in\R^n$,
lies entirely in $\lambda{\Omega}$ then by the definition of $g$
$$\int_{a+B_r}g(\Phi(x))dx = T(a).
$$

It remains to take care of the points $a$ such that $a+B_r$ does
not lie entirely  in $\lambda\Omega$. The dimension of the
boundary of $\Omega$ is less than or equal to  $n-\delta$ for some
$\delta>0$. Hence the number of balls of radius $r$ that are
required to cover the $r$-neighborhood of the boundary of
$\lambda\Omega$ (i.e. the set of all  points whose distance to the
boundary is less than $r$) is bounded by
$c\lambda^{n-\delta}$ for some constant $c$. The number of isolated
 points of $V$ in each
of such balls is bounded by the number of cells in our cell
decomposition. Therefore, the result of adding up $T$ over the
isolated points of $V$, whose distance to the boundary of
$\lambda\Omega$ is less than $r$, is bounded by
$c_0\lambda^{n-\delta}$ for some constant $c_0$. When we divide
the result by the volume of $\lambda\Omega$, which is proportional
to $\lambda^n$, the quotient approaches 0 as $\lambda$ approaches
infinity. Therefore, we can disregard the points that are close to
the boundary of $\lambda\Omega$ while computing the mean value
$M_{\Omega}$. We have proved that
$$\lim_{\lambda\to\infty}\frac{1}{\Vol(\lambda\Omega)}
\int_{\lambda\Omega}g(\Phi(x))dx $$ is the contribution of the
cell $C$ to the mean value $M_{\Omega}$.

We now work on the rewriting of the right hand side part of
\re{weyllaw}.

\begin{eqnarray}
\int_{\T^N}g(\varphi)d\varphi&=&\int_Q g(y,x)dydx = \int_{Q_y}
\left(\int_{Q_x} g(y,x) dx\right)dy\nonumber\\ &=& \int_{D} \w
T(y,H(y))dy=\int_C\w T(\varphi)p^*
d\varphi(\dots,A^1,\dots,A^n)=\int_C \w
T(\varphi)p^*\omega.\nonumber
\end{eqnarray}

Next we note that cells of dimension less than $N-n$ do not
contribute to the mean value. Indeed, for each of such cells we construct
a counting function $g$ and repeat the above argument.
It follows that the contribution of each of such cells is 
$$\int_{\T^N}g(\varphi)d\varphi,$$
which is equal to zero since $g$ is zero on a set of full
measure in $\T^N$.

The union of all 
$(i_1,\dots,i_{N-n},0,\dots ,0)$-cells that do not belong to the
closure of a cell of higher dimension is the $(N-n)$-dimensional
manifold $M_{N-n}(\w V)$, up to possibly a few cells of smaller
dimension. Those cells of smaller dimension do not contribute to
the mean value. We conclude that the mean value is equal to the
$\w T$-weighted $\cal O$-transversal volume of $\w V$.
\end{pf}

\section {Appendix A: Weyl's equidistribution law for dense
multidimensional orbits on the real torus.}

Classical Weyl's equidistribution law \cite{W3} states that
the time-average of a Riemann integrable function along a
one-dimensional dense orbit on a real torus coincides with the
space average. Here we state a version of Weyl's law for
multidimensional orbits.

Let $x\in\R^n$, $\varphi\in\T^N=\R^N/\Z^N$. A linear map
$\Phi:\R^n\to\R^N$ defines an action of $\R^n$ on $\T^N$ by
$\varphi\mapsto(\varphi+\Phi(x))\,\mod\,\Z^N$. Let $\cal O$ be an orbit
of this action.

\begin{Th}\label{T:Weyl}
Let $f$ be a  Riemann
integrable function on the torus $\T^N$, $\cal O$ an
orbit on $\T^N$ through an arbitrary point $\varphi_0\in\T^N$
such that there are no nonzero integral
vectors orthogonal to $\cal O$, and $\Omega$ a subset of
$\R^n$ with 
nonzero volume.

Then
$$\frac{1}{\Vol(\lambda\Omega)}\int_{\lambda\Omega}f(\varphi_0+\Phi(x))
dx\rightarrow \int_{\T^N} f(\varphi) d\varphi,\quad\text{as}\quad\lambda\rightarrow\infty.
$$
\end{Th}

The proof of this theorem in \cite{So} is a direct generalization of Weyl's
original argument for one-dimensional orbits.

\begin{Remark}\label{R:dense} This theorem implies that an
orbit is dense if and only if there are no nonzero integral
vectors orthogonal to the plane of the orbit.
\end{Remark}

\section{Appendix B: Cell-decomposition theorem}{\nonumber}
The cell-decomposition theorem for subanalytic sets
states that a subanalytic set can
be partitioned into finitely many cells, which are subanalytic
subsets of especially simple form. We will follow the formulation
of the cell-decomposition theorem as presented in \cite{vdD}.

\begin{Def}
Let $(i_1,\dots,i_n)$ be a sequence of zeros and ones of length
$n$. An $(i_1,\dots,i_n)$-{\it cell} is a subset of $\R^n$
obtained by induction on $n$ as follows:
\begin{itemize}
\item[(i)] a (0)-cell is a point in $\R$, a (1)-cell is an interval
$(a,b)\in\R$;
\item[(ii)] suppose $(i_1,\dots,i_n)$-cells are already defined; then
an $(i_1,\dots,i_n,0)$-cell is the graph $\Gamma(f)$ of a
$C^1$ function $f$ on an $(i_1,\dots,i_n)$-cell provided that $\Gamma(f)\subset\R^{n+1}$
is subanalytic in $\RP^{n+1}$;
further, an $(i_1,\dots,i_n,1)$-cell is a
set
$$(f,g)_X:=\{ (x,r)\in X\times\R: f(x)<r<g(x)\},
$$
where $X$ is an $(i_1,\dots,i_n)$-cell, functions $f,g$ are
$C^1$ functions on $X$ whose graphs $\Gamma(f),\Gamma(g)\subset\R^{n+1}$ are subanalytic in
$\RP^{n+1}$, and $f<g$ on $X$. The constant functions $f(x)=-\infty$,
$g(x)=+\infty$ are also allowed.
\end{itemize}
\end{Def}

\begin{Def} A {\it decomposition} of $\R^n$ is a special kind
of partition of $\R^n$ into finitely many cells. The definition is
by induction on $n$:
\begin{itemize}
\item[(i)] a decomposition of $\R$ is a collection
$$\left\{(-\infty,a_1),(a_1,a_2),\dots,(a_k,+\infty),\{a_1\},\dots,\{a_k\}\right\}
$$
where $a_1<\dots<a_k$ are points in $\R$.
\item[(ii)]a decomposition of $\R^{n+1}$ is a finite partition
into cells such that the set of projections of these cells to the
first $n$ coordinates forms a decomposition of $\R^n$.
\end{itemize}
\end{Def}

\begin{Th}\label{T:celldecomposition}{\sc (The cell decomposition theorem)}
Let $X$ be a subset of $\R^n$ subanalytic in $\RP^n$. Let $f$
be a function on $X$ whose graph in $\R^{n+1}$ is a subanalytic
subset of $\RP^{n+1}$. Then there exists a (finite) decomposition
$\cal D$ of $\R^n$ such that X is a union of cells in $\mathcal D$
and the restriction  $f|_B:B\to\R$ to each cell $B\in\cal D$ is
$C^1$.
\end{Th}

\end{document}